\documentclass[12pt]{article}
\usepackage[english]{babel}
\addto\captionsenglish{}
\usepackage[ansinew,latin1]{inputenc}
\usepackage[dvips]{graphicx}
\usepackage{amssymb,amsmath,amsfonts,amsthm}
\usepackage{latexsym,eucal,bbm,mathrsfs,xspace,nicefrac}
\usepackage[T1]{fontenc}

\usepackage{geometry}
\newcommand{\N}{\mathbbm{N}}
\newcommand{\Z}{\mathbbm{Z}}
\newcommand{\Q}{\mathbbm{Q}}
\newcommand{\R}{\mathbbm{R}}

\renewcommand{\phi}{\varphi}

\newcommand{\lm}{\Lambda M}

\newcommand{\refprop}[1]{Proposition~\ref{#1}}
\newcommand{\reflemma}[1]{Lemma~\ref{#1}}
\newcommand{\refcor}[1]{Corollary~\ref{#1}}

\DeclareMathOperator{\E}{E}
\DeclareMathOperator{\rank}{rank}

\DeclareMathOperator{\ind}{Index}
\DeclareMathOperator{\Ker}{Ker}

\DeclareMathOperator{\Hess}{Hess}
\DeclareMathOperator{\Fix}{Fix}
\DeclareMathOperator{\id}{id}

\DeclareMathOperator{\ex}{e}

\theoremstyle{plain}
\newtheorem{thm}{Theorem}
\newtheorem{thmi}{Theorem}
\newtheorem{prop}{Proposition}[section]
\newtheorem{lemma}{Lemma}[section]
\newtheorem{cor}{Corollary}[section]
\theoremstyle{definition}

\theoremstyle{remark}

\title{Three Dimensional Manifolds All of Whose
  Geodesics Are Closed} 

\author{John Olsen}

\begin{document}

\bibliographystyle{alpha}

\maketitle

\section{Introduction}

About $30$ years ago Berger conjectured that on a simply connected
manifold all of whose geodesics are closed, all geodesics have the
same least period. In addition to the spheres and projective spaces
with the standard metrics, the so-called Zoll metrics on $S^n$ have
this property as well; see \cite[Corollary 4.16]{Be;mgc}. The weaker
statement that there exists a common period is a special case of a
theorem due to Wadsley; see \cite[Theorem 7.12]{Be;mgc}. The lens
spaces with the canonical metrics show that simply connectedness is
necessary. On $S^{2n+1}/\Z_k$, $k > 2$ all geodesics are closed with
common period $2\pi$, but there exist geodesics of smaller period.

Bott and Samelson studied the topology of such manifolds and showed
that they must have the same cohomology ring as a compact rank one
symmetric space. In 1982 Gromoll and Grove proved the Berger
Conjecture for metrics on $S^2$, \cite[Theorem 1]{GG;S2}.

One possible way to attack the Berger Conjecture is to assume the
existence of shorter geodesics and then use Morse Theory on the free
loop space to derive a contradiction to the known (equivariant)
cohomology of the free loop space. In this paper we present some
results on the Morse theory on the free loop space of $S^3$ for
metrics all of whose geodesics are closed. In \cite{Z;FLSGSS} it is
shown that the energy function is perfect for the standard metrics on
all compact rank one symmetric spaces. The first theorem says that
with respect to $S^1$-equivariant cohomology, the same is true for all
metrics on $S^3$ all of whose geodesics are closed.

\begin{thmi}
  Let $g$ be a metric on $S^3$ all of whose geodesics are closed. Then
  relative to the point curves the energy function $\E$ is perfect
  with respect to rational $S^1$-equivariant cohomology.
\end{thmi}

The fact that the energy function is a Morse-Bott function was
observed by Wilking; see \cite[Proof of Step 3]{W;Index}. This allows
us to do Morse Theory on $\Lambda S^3$. The next theorem allows us to
use arbitrary coefficients for the homology.

\begin{thmi}
  Let $g$ be a metric on $S^3$ all of whose geodesics are closed. Then
  the negative bundles over the critical manifolds are orientable both
  as ordinary and $S^1$-equivariant vector bundles.
\end{thmi}

We have the following structure result for the critical manifolds of
the energy function.

\begin{thmi}
  Let $g$ be a metric on $S^3$ all of whose geodesics are closed. A
  critical manifold in $\Lambda S^3$ is either diffeomorphic to the
  unit tangent bundle of $S^3$ or has the integral cohomology ring of
  a three dimensional lens space $S^3/\Z_{2k}$.
\end{thmi}

Together with the Bott Iteration Formula for the index of the iterates
of a closed geodesic, the above results have strong implications for
the geometry of the manifold; see Section~\ref{conseq}. 

The first part of this work was advised by Burkhard Wilking at the
University of M\"unster. I would like to thank him for introducing me
to the problem and for helpful advice how to attack it. The second
part of this work was advised by Wolfgang Ziller at the University of
Pennsylvania. I would like to thank him for generous and helpful
advice. I would also like to thank Karsten Grove and Marcel B\"okstedt
for many helpful conversations.

\section{Preliminaries}\label{prelim}

We will review some basic notions from Morse theory on the free loop
space of a Riemannian manifold. The standard reference is
\cite{K;LCG}; see also \cite[Section 3]{BO;SS} and \cite{H;Eqmorse}
for a short introduction.

Let $M$ be a Riemannian manifold and let $\lm = W^{1,2}(S^1,M)$ be the
free loop space of $M$. The free loop space $\lm$ can be given the
structure of a smooth Hilbert manifold which makes $\E$ into a smooth
function; see \cite[Theorem 1.2.9]{K;LCG}. Note that there is a
natural inclusion $W^{1,2}(S^1,M) \to C^0(S^1, M)$ and a fundamental
theorem states that this map is a homotopy equivalence; see
\cite[Theorem 1.2.10]{K;LCG}. The tangent space at a point $c \in
\lm$, $T_c\lm$, consists of all vector fields along $c$ of class
$W^{1,2}$, and is a real Hilbert space with inner product given by
\begin{align*}
  \langle \langle X,Y \rangle \rangle_1 = \int_0^{2\pi} \langle X(t),
  Y(t)\rangle + \langle \frac{DX}{dt},\frac{DY}{dt} \rangle dt,
\end{align*}
where $\frac{D}{dt}$ is the covariant derivative along $c$
induced by the Levi-Civita connection on $M$. Let the energy function
$\E \colon \lm \to \R$ be given by $\E(c) =
\int_0^{2\pi}|\dot c(t)|^2dt$. It follows from the first variation
formula that the critical points of $\E$ are the closed geodesics on
$M$.

Let $N$ be a critical manifold of $\E$ in $\lm$ and let $c \in N$ be a
critical point. We use the convention that the critical sets are
maximal and connected throughout. By definition $N$ satisfies the Bott
nondegeneracy condition if $T_cN = \Ker(\Hess_c(\E))$. If all critical
manifolds are nondegenerate in this sense, the energy function is
called a Morse-Bott function. What this says geometrically is that the
dimension of the space of periodic Jacobi vector fields along the
geodesic $c$ is equal to the dimension of the critical manifold.

Now assume that $N_j$, $j = 1, \ldots, l$, are the critical manifolds
of energy $a$, that they satisfy the Bott nondegeneracy condition, and
that there are no other critical values in the interval $[a-\epsilon,
a +\epsilon]$, $\epsilon > 0$. The metric on $T_c \lm$ induces a
splitting of the normal bundle, $\xi$, of $N_j$ in $\lm$ into a
positive and a negative bundle, $\xi = \xi^+ \oplus \xi^-$, such that
the Hessian of the energy function is positive definite on $\xi^+$ and
negative definite on $\xi^-$. Furthermore $\lambda(N_j) = \rank \xi^-$
is called the index of $N_j$ and is finite. Note that the index is
constant on each critical manifold, since it is connected. We denote
by $\Lambda^aM$ the set $\E^{-1}([0,a]) \subseteq \lm$. Let
$D(\xi^-(N_j)) = D(N_j)$ and $S(\xi^-(N_j))=S(N_j)$ be the disc,
respectively sphere, bundle of the negative bundle $\xi^-$ over
$N_j$. Using the gradient flow, one shows that there exists a homotopy
equivalence $\Lambda^{a+ \epsilon}M \simeq \Lambda^{a- \epsilon}M
\cup_f \cup_{j=1}^l D(N_j)$ for some gluing maps $f_j \colon S(N_j)
\to \Lambda^{a- \epsilon}M$; see \cite[Theorem 2.4.10]{K;LCG}.

If the negative bundle over $N_j$ is orientable for all $j$, excision
and the Thom isomorphism yield
\[
H^i(\Lambda^{a+ \epsilon}M, \Lambda^{a- \epsilon}M; R) \cong
\bigoplus_{j=1}^l H^i(D(N_j), S(N_j); R) \cong \bigoplus_{j=1}^l
H^{i-\lambda(N_j)}(N_j; R)
\]
for any coefficient ring $R$. If $N$ is not orientable the isomorphism
holds with $\Z_2$ coefficients. The cohomology of $\Lambda^{a+
  \epsilon}M$ is now determined by the cohomology of $\Lambda^{a-
  \epsilon}M$ and $N_j$ by the long exact cohomology sequence for the
pair $(\Lambda^{a+ \epsilon}M, \Lambda^{a- \epsilon}M)$. This way one
can in principle inductively calculate the cohomology of $\lm$ from
the cohomology of the critical manifolds. If the map
$\Lambda^{a+\epsilon}M \to (\Lambda^{a+\epsilon}M, \Lambda^{a-
  \epsilon}M)$ induces an injective map $H^i(\Lambda^{a+ \epsilon}M,
\Lambda^{a-\epsilon}M; R) \to H^i(\Lambda^{a+ \epsilon}M; R)$ for all
$i$, we say that all relative classes can be completed to absolute
classes. This is equivalent to all the boundary maps in the long exact
cohomology sequence for the pair $(\Lambda^{a+ \epsilon}M, \Lambda^{a-
  \epsilon}M)$ being zero. If this holds for all $i$ and all critical
values $a$, we say that $\E$ is perfect. If $R$ is a field,
perfectness implies that $H^j(\lm; R) = \bigoplus_j H^{j -
  \lambda(N_j)}(N_j; R)$, where the sum is over all critical
manifolds.

For a topological group $G$ we let $EG$ be a contractible topological
space on which $G$ acts freely and let $EG/G=BG$ be the classifying
space of $G$. For a $G$-space $X$ the quotient space $(X \times EG)/G
= X \times_G EG = X_G$ is called the Borel construction. The
$G$-equivariant cohomology of $X$ is defined to be the usual
cohomology of the Borel construction $X_G$. We have a fibration $X \to
X \times_G EG \to BG$ and, if the $G$-action on $X$ is free, a
fibration $X \times_G EG \to X/G$ with fiber $EG$, i.e. the map is a
(weak) homotopy equivalence. If $G$ acts on a manifold $X$ with finite
isotropy groups, the map $X \times_G EG \to X/G$ is a rational
homotopy equivalence and the cohomology of $X/G$ is concentrated in
finitely many degrees when $X$ is finite dimensional and compact.

For $G$ a group acting on $\lm$ we consider $G$-equivariant Morse
theory. Let $\xi \to X$ be a vector bundle where $G$ acts on $\xi$
such that $p\colon \xi \to X$ is equivariant and the action is linear
on the fibers. For such vector bundles we define the $G$-vector bundle
by $\xi_G= \xi \times_G EG \to X_G$. The $G$-vector bundle $\xi_G$ is
orientable if and only if $\xi$ is orientable and $G$ acts orientation
preserving on the fibers. Thus if $G$ is connected $\xi_G$ is
orientable if and only if $\xi$ is orientable. For an oriented rank
$k$ $G$-vector bundle $\xi_G$ over $X$ there is a $G$-equivariant Thom
isomorphism, $H^*_G(D((X)), S((X)); R) \cong H^{*-k}_G(X; R)$.

The action of $O(2)$ on $S^1$ induces an action of $O(2)$ on the free
loop space via reparametrization. Since the energy function is
invariant under the action of $O(2)$ and $O(2)$ acts by isometries
with respect to the inner product used to define the negative bundles,
the negative bundles are $O(2)$-bundles in this sense. For any group
$G \subseteq O(2)$ we have similarly to the above
\[
H^i_G(\Lambda^{a+ \epsilon}M, \Lambda^{a- \epsilon}M; R) \cong
\bigoplus_{j=1}^l H^i_G(D(N_j),S(N_j); R) \cong \bigoplus_{j=1}^l
H^{i-\lambda(N_j)}_G(N_j; R)
\]
where as above $D(N_j)$ and $S(N_j)$ are the disk and sphere bundle of
the negative bundles $\xi^-$ over the critical manifold $N_j$ $j = 1,
\ldots, l$ and where $a$ is the energy of the critical manifolds
$N_j$. The second isomorphism holds for any coefficient ring $R$ if
the negative bundle is oriented and with $\Z_2$ coefficients if it is
nonorientable. The $G$-equivariant cohomology of $\Lambda^{a+
  \epsilon}M$ is determined by the $G$-equivariant cohomology of
$\Lambda^{a- \epsilon}M$ and $N_j$ by the long exact $G$-equivariant
cohomology sequence for the pair $(\Lambda^{a+ \epsilon}M, \Lambda^{a-
  \epsilon}M)$. In principle this allows us to inductively calculate
the $G$-equivariant cohomology of $\lm$ from the $G$-equivariant
cohomology of the critical manifolds. Again we say that $\E$ is
perfect with respect to $G$-equivariant cohomology if the map
$\Lambda^{a+\epsilon}M \to (\Lambda^{a+\epsilon}M, \Lambda^{a-
  \epsilon}M)$ induces an injective map $H_G^i(\Lambda^{a+ \epsilon}M,
\Lambda^{a-\epsilon}M; R) \to H_G^i(\Lambda^{a+ \epsilon}M; R)$ for
all $i$ and for all critical values $a$. If $R$ is a field,
perfectness implies that $H_G^j(\lm; R) = \bigoplus_j H_G^{j -
  \lambda(N_j)}(N_j; R)$, where the sum is over all critical
manifolds. See \cite{H;Eqmorse} for an introduction to equivariant
Morse theory on $\lm$.

\section{Rational $S^1$-equivariant Perfectness of the Energy
  Function}\label{perfectness}

In the case of the canonical metric $g_0$ on $S^3$ we let $B_k \cong
T^1S^3$ be the manifold of $k$-times iterated geodesics. The energy
function is a Morse-Bott function for the metric $g_0$. The manifolds
$B_k$, $k \in \N_{>0}$ and the point curves $S^3$ are the only
critical manifolds for $\E$ and the induced action of $S^1/\Z_k$ on
$B_k$ is free. The $S^1$-equivariant cohomology of $T^1S^3$ can thus
be calculated from the Gysin sequence for the bundle $S^1 \to T^1S^3
\to T^1S^3/S^1$ and is given by
\[
H^i_{S^1}(T^1S^3; \Q) =
\begin{cases}
  \Q, \quad & i=0,4,\\
  \Q^2, \quad & i =2,\\
  0, \quad & \text{otherwise}.
\end{cases}
\]
Since by \cite{Z;FLSGSS} the indices of the critical manifolds $B_k$
in $\Lambda S^3$ are $2(2k-1)$, $k\in \N_{>0}$, we see that for all
$k$ the rational $S^1$-equivariant cohomology of $B_k$ is concentrated
in even degrees. Hence by the Lacunary Principle the energy function
is perfect relative to $S^3$ for rational $S^1$-equivariant
cohomology.

\begin{prop}\label{s1eqcohom}
  The rational $S^1$-equivariant cohomology of $\Lambda S^3$ relative
  to $S^3$ is given by
  \[
  H^i_{S^1}(\Lambda S^3,S^3; \Q) =
  \begin{cases}
    \Q, \quad &i = 2\\
    \Q^2, \quad &\text{$i=2k$, $k \in \N$, $k>1$}\\
    0, \quad & \text{otherwise}.
  \end{cases}
  \]
\end{prop}

More generally, Hingston calculated the $S^1$-equivariant cohomology
of the free loop space of any compact rank one symmetric space; see
\cite[Section 4.2]{H;Eqmorse}.

Let $g$ be a metric on $S^3$ all of whose geodesics are closed and
normalized so that $2\pi$ is the least common period. We assume that
there exist exceptional geodesics on $(S^3, g)$ of period $2\pi/n$ for
$n > 1$. Since all geodesics are closed with common period $2\pi$ the
geodesic flow defines an effective and orientation preserving action
of $\R/2\pi \Z = S^1$ on $T^1S^3$. Choose a metric on $T^1S^3$ such
that the action of $S^1$ becomes isometric. The full unit tangent
bundle corresponds to geodesics of period $2\pi$ and the closed
geodesics of length $(k/n)2\pi$ can be identified with the fixed point
set of the element $e^{2\pi i /n} \in S^1$ in $T^1S^3$, i.e. the fixed
point set of $\Z_n \subseteq S^1$. The fixed point set of $\Z_n
\subseteq S^1$ again has an effective and orientation preserving
action of $S^1/\Z_n = S^1$ since $S^1$ is abelian. Since the critical
sets of the energy function can be identified with the fixed point
sets of some $g \in S^1$ and the metric is chosen so that the action
is isometric, it follows that the critical sets are compact, totally
geodesic submanifolds of $T^1S^3$.

We recall an observation made by Wilking, \cite[Proof of Step
3]{W;Index}: If all geodesics are closed, the energy function is a
Morse-Bott function. To see this, one notes that a critical manifold
$N \subseteq T^1S^3$ of geodesics of length $2\pi/n$ is a connected
component of the fixed point set of an element $g = \ex^{2\pi i/n} \in
S^1$. The dimension of the critical manifold is equal to the
multiplicity of the eigenvalue $1$ of the map $g_{*v}$ at a fixed
point $v$. Since the differential of the geodesic flow is the
Poincar\'{e} map, the multiplicity of the eigenvalue $1$ at $v$ is
equal to the dimension of the vector space of $2\pi/n$-periodic Jacobi
fields along the geodesic $c (t)= \exp(tv)$, $t \in [0,2\pi/n]$ ($\dot
c$ is considered as a periodic Jacobi field as well). Since the null
space of $\Hess_c(\E)$ consists of $2\pi/n$-periodic Jacobi vector
fields, we see that the kernel of $\Hess_c(\E)$ is equal to the
tangent space of the critical manifolds.

We know that the geodesic flow acts orientation preserving on the
fixed point sets, but we need to show that the fixed point sets are
indeed orientable. This and the fact that the fixed point sets have
even codimension follow from the following lemma.

\begin{lemma}\label{fixorient}
  The fixed point sets $\Fix(\Z_n) \subseteq T^1S^3$ have even
  codimension and are orientable.
\end{lemma}

\begin{proof}
  We recall a few facts about the Poincar\'{e} map; see \cite[page
  216]{BTZ;closedgeo}. Let $c$ be a closed geodesic with $c'(0) =v
  \neq 0$. Let $\phi_t$ denote the geodesic flow and note that a
  closed geodesic with $v = c'(0)$ corresponds to a periodic orbit
  $\phi_tv$. The flow $\phi_t$ maps the set $T_rS^3 = \{v \in TS^3
  \mid |v|=r \}$ into itself. The Poincar\'{e} map $\mathcal{P}$ of
  $c$ is the return map of a local hypersurface $N \subset T_rS^3$
  transversal to $v$. The linearized Poincar\'{e} map is up to
  conjugacy independent of $N$ and is given by $P = D_v
  \mathcal{P}$. One can choose $N$ such that $T_vN = V \oplus V$,
  where $V = v^{\perp} \subseteq T_pS^3$, where $p$ is the footpoint
  of $v$. On $V \oplus V$ there is a natural symplectic structure
  which is preserved by $P$. As above we identify the critical
  manifolds with the fixed point sets of $\Z_n \subset S^1$ and notice
  that since the critical manifolds are nondegenerate, the dimension
  of the fixed point sets is equal to the multiplicity of the
  eigenvalue $1$ of $P$ plus one. By \cite[Proposition 3.2.1]{K;LCG}
  we know that the multiplicity of $1$ as an eigenvalue of $P$ is
  even. Hence the dimension of $\Fix(\Z_n) \subseteq T^1S^3$ is odd
  and the codimension is even.
  
  If the fixed point sets are one or five dimensional, they are
  clearly orientable, so assume $\dim \Fix(\Z_n) = 3$. In this case
  the normal form of $P$ has a $2\times 2$ identity block and a $2
  \times 2$ block which is a rotation, possibly by an angle
  $\pi$. Denote the two dimensional subspaces by $V_{\id}$ and
  $V_{\text{rot}}$. Using the symplectic normal form for $P$ we know
  that $V_{\id}$ and $V_{\text{rot}}$ are orthogonal with respect to
  the symplectic form and that the restriction of the symplectic form
  to each subspace is nondegenerate; see \cite[page
  222]{BTZ;closedgeo}. If we consider the normal bundle $\nu
  \Fix(\Z_n)$ with fiber $\nu_v \Fix(\Z_n)$ for $v \in \Fix(\Z_n)$, we
  have $\nu_v \Fix(\Z_n) = V_{\text{rot}}$. Since $V_{\text{rot}}$ is
  a symplectic subspace, it carries a canonical orientation. This
  gives a canonical orientation on each fiber of the normal bundle,
  which shows that the normal bundle is orientable and hence that
  $\Fix(\Z_n)$ is orientable.
\end{proof}

The one dimensional critical manifolds are diffeomorphic to
circles. If $c$ is a one dimensional critical manifold, the geodesic
$\overline c(t) =c(-t)$ is a second critical manifold of the same
index. If the critical manifold is five dimensional it is the full
unit tangent bundle. If the exceptional critical manifold that
consists of prime closed geodesics is five dimensional then all
geodesics are closed with period $2\pi/n$ contradicting the
assumption.

We now state and prove the main result of the thesis.

\begin{thm}\label{perfect}
  Let $g$ be a metric on $S^3$ all of whose geodesics are closed. Then
  relative to the point curves the energy function $\E$ is perfect
  with respect to rational $S^1$-equivariant cohomology.
\end{thm}

\begin{proof}
  The proof takes up several pages, and we first give a short
  outline. We will use the Index Parity Theorem repeatedly; see
  \cite[Theorem 3]{W;Index}. The Index Parity Theorem states that for
  an oriented Riemannian manifold $M^n$ all of whose geodesics are
  closed, the index of a geodesic in the free loop space is even if
  $n$ is odd, and odd if $n$ is even; in particular, in our case it
  states that all indices are even. The idea is to show that the
  contributions $H^i_{S^1}(D(N), S(N); \Q)$ from a critical manifold
  $N$ to the $S^1$-equivariant cohomology of $\Lambda S^3$ occur in
  even degrees only. If the negative bundle over $N$ is orientable
  this is by the Thom Isomorphism equivalent to showing that
  $H^i_{S^1}(N; \Q) = 0$ for $i$ odd. Perfectness of the energy
  function then follows from the Lacunary Principle, since the
  critical manifolds all have even index.

  We first show that the negative bundles over the one and five
  dimensional critical manifolds are oriented and that the critical
  manifolds only have $S^1$-equivariant cohomology in even degrees. If
  the critical manifold is three dimensional we first use Smith Theory
  to show that the $S^1$-Borel construction is rationally homotopy
  equivalent to $S^2$.  If the negative bundle is oriented the
  contributions occur in even degrees only. If the negative bundle is
  nonorientable we use Morse Theory and a covering space argument to
  show that the critical manifold does not contribute to the
  $S^1$-equivariant cohomology.

  We begin the proof by considering the five dimensional case. The
  five dimensional critical manifold is diffeomorphic to the unit
  tangent bundle. It is clear that the negative bundles over the unit
  tangent bundle are oriented, since the unit tangent bundle is simply
  connected. The Borel construction $T^1S^3 \times_{S^1} ES^1$ is
  rationally homotopy equivalent to $T^1S^3/S^1$, since the action has
  finite isotropy groups, and thus the possible degrees where $T^1S^3
  \times_{S^1}ES^1$ has nonzero rational cohomology is zero through
  four. Using these facts and the Gysin sequence for the bundle $S^1
  \to T^1S^3 \times ES^1 \to T^1S^3 \times_{S^1} ES^1$ we see that
  $H^*_{S^1}(T^1S^3; \Q) = H^*(S^2 \times S^2; \Q)$ and hence has
  nonzero classes in even degrees only.

  Next, we show that the negative bundles over the one dimensional
  critical manifolds are oriented and thus, since the $S^1$-Borel
  construction of the one dimensional critical manifolds has rational
  cohomology as a point, the contributions occur in even degrees only.

  \begin{prop}\label{orient}
    The negative bundles over the one dimensional critical manifolds
    are ($S^1$-equivariantly) orientable.
  \end{prop}

  \begin{proof}
    The negative bundle over the critical manifold consisting of a
    prime closed geodesics $c$ is orientable, since we can define an
    orientation in one fiber of the negative bundle and use the free
    $S^1/\Z_n$-action to define an orientation in the other fibers.

    Let $\xi_n$ be the negative bundle over the $n$-times iterated
    geodesic. The representation of $\Z_n$ on the fibers of $\xi_n$ is
    presented in \cite[Proposition 4.1.5]{K;LCG}. The representation
    of $\Z_n$ is the identity on a subspace of dimension $\ind(c)$
    (which corresponds to the image of the bundle $\xi_1$ under the
    $n$-times iteration map) and is a sum of two dimensional real
    representations given by multiplication by $\ex^{\pm 2\pi ip/n}$
    on a vector space of even dimension. If $n$ is even, $\Z_n$ also
    acts as $-\id$ on a subspace of dimension $\ind(c^2) - \ind(c)$.

    We know by \cite[Proposition 4.1.5]{K;LCG} that the dimension of
    the subspace on which a generator $T$ of $\Z_n$ acts as $-\id$ is
    equal to $\ind(c^2) - \ind(c)$, which by the Index Parity Theorem
    is even. By \cite[Lemma 4.1.4]{K;LCG} the pair $(D^k/\Z_n,
    S^{k-1}/\Z_n)$ is orientable ($k = \ind(c^n)$), since the
    dimension of the subspace on which $T$ acts as $-\id$ is even
    dimensional. Pick an orientation in one fiber of the negative
    bundle invariant under the action of $\Z_n$. Use the $S^1/\Z_n$-
    action to define an orientation in any other fiber of the negative
    bundle. Since the orientation is chosen to be invariant under the
    action of $\Z_n$ the orientation on the negative bundle is
    well-defined.
  \end{proof} 

  We now treat the three dimensional case and first prove the
  following important fact.
 
  \begin{prop}\label{eqcofix}
    Assume that the critical manifold $N$ has dimension three. Then
    the $S^1$-Borel construction of $N$ is rationally homotopy
    equivalent to $S^2$.
  \end{prop}

  \begin{proof}
    Assume that $\Z_n$ acts trivially on $N$. We start by considering
    the quotient $N \to N/(S^1/\Z_n)$. The action of $S^1/\Z_n$ is
    effective and has isotropy at a closed geodesic which is an $l$th
    iterate, for some $l$, of a closed geodesic in a one dimensional
    critical manifold. The $S^1/\Z_n$-action normal to the
    $S^1/\Z_n$-orbit acts as a rotation by $\ex^{2\pi i/l}$. Hence
    $N/(S^1/\Z_n)$ is a two dimensional orbifold, since a neighborhood
    of an arbitrary point in $N/(S^1/\Z_n)$ is homeomorphic to $\R^2/
    \Z_l$. However, the quotient $\R^2/ \Z_l$ is homeomorphic to
    $\R^2$ so in particular $N/ (S^1/\Z_n)$ is homeomorphic to a
    surface. Since by \reflemma{fixorient} $\Fix(\Z_n)$ is orientable
    and the action of $S^1/\Z_n$ is orientation preserving, the
    quotient is an orientable surface of genus $g$.

    Using the Gysin sequence we can calculate $H^1(N; \Q)$ from the
    bundle $S^1 \to N \times ES^1 \to N \times_{S^1} ES^1$, since we
    know that $N/(S^1/\Z_n)$ is rationally homotopy equivalent to $N
    \times_{S^1} ES^1$. The Gysin sequence yields
    \[
    H^i(N; \Q) =
    \begin{cases}
      \Q \quad &i = 0,3\\
      \Q^{2g}, \quad & i = 1, 2 \text{ if } \chi \neq 0,\\
      \Q^{2g+1}, \quad &i = 1, 2 \text{ if } \chi = 0,
    \end{cases}
    \]
    where $\chi$ denotes the Euler class of the bundle.

    If $n$ is a prime, $p$, we can use Smith theory for the
    $\Z_n$-action on $T^1S^3$ to bound the sum of the Betti numbers of
    the fixed point sets, i.e. we know that the sum $\sum
    b_i(\Fix(\Z_p); \Z_p) \leq \sum b_i(T^1S^3; \Z_p) = 4$ for an
    arbitrary prime $p$; see \cite[Theorem 4.1]{Br;compgr}. Hence we
    have $\sum b_i(\Fix(\Z_p); \Q) \leq 4$, by the Universal
    Coefficient Theorem. If $n$ is not prime we choose a $p$ such that
    $p | n$ and such that $\Fix(\Z_p)$ is three dimensional, which is
    possible since $N$ is three dimensional and the only critical
    manifold of dimension five is $T^1S^3$. Then $N \subseteq
    \Fix(\Z_p)$ and since both $N$ and $\Fix(\Z_p)$ are closed three
    dimensional submanifolds of $T^1S^3$, $N$ equals a component of
    $\Fix(\Z_p)$. Hence we have $\sum b_i(N; \Q) \leq 4$, which
    implies that $g=0$.
  \end{proof}

  Thus, if the negative bundle is orientable we only have
  contributions in even degrees. We now treat the case where the
  negative bundle is nonorientable.

  \begin{prop}\label{nonor}
    Let $\xi^- \to N$ be a nonorientable negative bundle over a three
    dimensional critical manifold $N$. Then the cohomology groups
    $H^*_{S^1}(D(N), S(N); \Q)$ vanish.
  \end{prop}

  \begin{proof}
    Let $p\colon \tilde N \to N$ be the twofold cover such that the
    pull-back $p^* \xi^- = \tilde \xi^-$ is orientable. Let $\iota$ be
    the covering involution on $\tilde N$ which lifts to $\tilde
    \xi^-$. Note that the $S^1$-equivariant negative bundle
    $\xi^-_{S^1}\to N_{S^1}$ is also nonorientable. Lift the action of
    $S^1/\Z_n$, or possibly a twofold cover of $S^1/\Z_n$, on $\xi^-$
    to an action on $\tilde{\xi}^-$ such that $\iota$ becomes
    equivariant with respect to this action. Let $\tilde{\xi}^-_{S^1}
    \to \tilde N_{S^1}$ be the oriented twofold cover of $\xi_{S^1}
    \to N_{S^1}$. We know from \refprop{eqcofix} that $N_{S^1}$ is
    rationally homotopy equivalent to $S^2$. The action of $S^1/\Z_n$
    on $\tilde N$ also has finite isotropy groups, so by an argument
    similar to the one in the proof of \refprop{eqcofix} we see that
    $\tilde N /(S^1/\Z_n)$ is an oriented surface of genus $g$,
    $F_g$. The Borel construction $\tilde N_{S^1}$ is rationally
    homotopy equivalent to $F_g$.

    There are two cases to consider: $g=0$ and $g > 0$. We first treat
    the case $g= 0$. Since the bundle $\tilde \xi^-_{S^1}$ is
    orientable there exists a Thom class and this cohomology class
    changes sign under the action of $\iota^*$, since otherwise it
    would descend to give a Thom class for $\xi^-_{S^1}$ making it
    orientable. By general covering space theory we know that
    $H^*_{S^1}(D(\tilde N), S(\tilde N); \Q)^{\iota^*} =
    H^*_{S^1}(D(N), S(N); \Q)$ and $H^*_{S^1}(\tilde N; \Q)^{\iota^*}
    = H^*_{S^1}(N; \Q)$, where the superscript $\iota^*$ denotes the
    classes that are fixed under the action of $\iota^*$, and where as
    before, $D(N) = D(\xi(N))$ etc.

    By the Thom isomorphism we have $H^i_{S^1}(D(\tilde N), S(\tilde
    N); \Q) = H^{i - \lambda(N)}_{S^1}(\tilde N; \Q)$ and since $g=0$
    we have $H^i_{S^1}(\tilde N; \Q) =0$ for all $i \neq 0,2$, which
    means that we only have to see what happens to the degree zero and
    degree two cocycles. By the above, this implies that
    $H^i_{S^1}(D(N), S(N); \Q)= 0$ for $i \neq \lambda(N), \lambda(N)
    +2$. The action of $\iota^*$ on $H^0_{S^1}(\tilde N; \Q)$ is
    trivial and since the Thom class changes sign, the action of
    $\iota^*$ on $H^{\lambda(N)}_{S^1}(D(\tilde N),S(\tilde N); \Q)$
    has no fixed points. Hence we conclude that
    $H^{\lambda(N)}_{S^1}(D(N),S(N); \Q) =0$.

    Since both $N$ and $\tilde N$ are two spheres we have
    $H^2_{S^1}(\tilde N; \Q) = \Q$ and $H^2_{S^1}(N; \Q) = \Q$, but
    also that $H^2_{S^1}(\tilde N); \Q)^{\iota^*} = H^2_{S^1}(N; \Q) =
    \Q$. This implies that $\iota^*$ acts trivially on
    $H^2_{S^1}(\tilde N; \Q)$ and, since the Thom class changes sign
    under the action of $\iota^*$, as $-\id$ on
    $H^{\lambda(N)+2}_{S^1}(D(\tilde N), S(\tilde N); \Q)$. We
    conclude that $H^{\lambda(N)+2}_{S^1}(D(N), S(N); \Q) = 0$.

    The second case to consider is $g > 0$. In that case we will
    derive a contradiction. Let $2k$ be the minimal index of a
    critical manifold on which $\Z_n$ acts trivially and whose
    negative bundle is nonorientable. If there is more than one
    critical manifold of index $2k$ with nonorientable negative
    bundle, we repeat the argument for each manifold. If the action of
    $S^1/\Z_n$ had been free, $N_{S^1}$ would be homotopy equivalent
    to $S^2$ and the negative bundle would have been orientable. Hence
    the action of $S^1/\Z_n$ on $N$ has fixed points and since the
    action is effective, the fixed points are one dimensional critical
    manifolds.

    Since the map $\tilde N/S^1 \to N/S^1$ is a branched covering, the
    Riemann-Hurwitz formula implies that there are $2g +2$ branched
    points.  Exceptional orbits of the $S^1$-action on $\tilde N$ are
    circles, since the action is effective, and they project down to
    exceptional orbits (shorter geodesics) for the action on
    $N$. Branched points for the covering $\tilde N / S^1 \to N/S^1$
    correspond to orbits where the isotropy of the action on $N$ is
    bigger than for the action on $\tilde N$. Hence we see that the
    branched points come from exceptional orbits for the action on $N$
    and hence that there are at least $2g+2$ exceptional orbits in
    $N$.

    The exceptional orbits are shorter geodesics, and since by the
    Bott Iteration Formula $\ind (c^q) \geq \ind (c)$, these circles
    must all have index less than or equal to $2k$. By the following
    lemma we deduce that the index must be equal to $2k$. For later
    reference we state this lemma separately.

    \begin{lemma}\label{nonor-claim}
      There exists no one dimensional critical manifolds of index less
      than $2k$.
    \end{lemma}

    \begin{proof}
      Assume for contradiction that there exists a one dimensional
      critical manifold of index $2h < 2k$, $h \geq 2$. This critical
      manifold hence contributes a $\Q^2$ in degree $2h$. A three
      dimensional critical manifold contributes a $\Q$ in the degree
      equal to the index and a $\Q$ in degree equal to the index plus
      $2$, since the negative bundle is orientable (hence a total of
      two classes). The unit tangent bundle contributes a total of $4$
      classes. There are no cancellations in degree less than $2k$,
      since the contributions all have even degrees. The total number
      of classes needed to "fill the gaps" (from $H_{S^1}^2(\Lambda
      S^3, S^3; \Q)$ to $H_{S^1}^{2h-2}(\Lambda S^3, S^3; \Q)$) in the
      cohomology is $2h-3$, since $H^2_{S^1}(\Lambda S^3, S^3; \Q) =
      \Q$. Notice that as $2h < 2k$ no cancellations are possible,
      since all contributions occur in even degrees. This implies that
      no three dimensional critical manifold can have index $2h-2$ and
      the unit tangent bundle cannot have index $2h-2$ or $2h-4$,
      since otherwise $\dim H_{S^1}^{2h}(\Lambda S^3, S^3; \Q) \geq
      3$, since we have contributions from the two circles and the
      three or five dimensional critical manifold. Hence we have to
      fill an odd number of holes with contributions that only come in
      pairs or quartets. This is clearly impossible.

      If $h =1$ we get a contribution of $\Q^2$ in degree $2$, which
      cannot cancel out, since the index two critical manifold only
      contributes in even degrees. This contradicts the fact that
      $H^2_{S^1}(\Lambda S^3, S^3; \Q) = \Q$.
    \end{proof}

    By a similar argument, we see that there is at most one three or
    five dimensional critical manifold for each even index less than
    $2k$. If the unit tangent bundle has index $2j$, there is no
    critical manifold of index $2j + 2$, since $H_{S^1}^2(T^1S^3; \Q)
    = \Q^2$. 

    Now we derive the contradiction. This is done by considering the
    possible contributions to $H^{2k}_{S^1}(\Lambda S^3, S^3;
    \Q)$. The three dimensional critical manifold of index $2k-2$
    contributes a $\Q$ in degree $2k$, a five dimensional critical
    manifold of index $2k-4$ contributes a $\Q$ in degree $2k$, and
    similarly a five dimensional critical manifold of index $2k-2$
    contributes a $\Q^2$ in degree $2k$. Hence we get a contribution
    to $H^{2k}_{S^1}(\Lambda S^3, S^3; \Q)$ of at least a $\Q$ from
    the critical manifold of index $2k-2$ or $2k-4$ and a $\Q^{2g+2}$
    from the exceptional orbits. By the Thom isomorphism we have
    $H^{2k +1}_{S^1}(D(\tilde N), S(\tilde N); \Q) = H^1_{S^1}(\tilde
    N; \Q)$ and furthermore we have by covering space theory that
    $H^{2k+1}_{S^1}(D(\tilde N), S(\tilde N); \Q)^{\iota^*} = H^{2k
      +1}_{S^1}(D(N), S(N); \Q)$, which implies that $H^1_{S^1}(\tilde
    N; \Q)^{\iota^*} = H^{2k+1}_{S^1}(D(N), S(N); \Q)$. Since
    $H^1_{S^1}(\tilde N; \Q) = \Q^{2g}$, we see that the maximal
    contribution to $H^{2k+1}_{S^1}(\Lambda S^3, S^3; \Q)$ is
    $\Q^{2g}$, which for example happens if $\iota^*$ acts as $-\id$
    on $\Q^{2g}$. This yields a contradiction since we now have $\dim
    H^{2k}_{S^1}(\Lambda S^3, S^3; \Q) \geq 3$.

    If there exists another critical manifold of index $2k$ with
    nonorientable negative bundle we repeat the argument above. The
    exceptional geodesics that contribute in degree $2k$ are distinct
    from the ones in other critical manifolds since the iterates of
    the shorter geodesics lie in different connected critical
    manifolds.
  \end{proof}

  This finishes the proof that $\E$ is perfect.
\end{proof}

A consequence of the perfectness of the energy function is the
following general fact.

\begin{cor}\label{conn}
  The minimal index of a critical manifold consisting of nonconstant
  geodesics is two. For every number $2k$ there exists at most one
  connected critical manifold of index $2k$ for every $k \geq 1$.
\end{cor}

\begin{proof}
  As a global minimum of $\E$, $S^3$ has index zero. If one of the
  critical manifolds $N$ consisting of nonconstant geodesics has
  index zero, $H_{S^1}^0(\Lambda S^3; \Q)$ would be at least two
  dimensional, which is not the case. If the minimal index, $\ind(N) =
  2i$, $i > 1$, then $H^2(\Lambda S^3; \Q) = \cdots = H^{2i-1}
  (\Lambda S^3;\Q) =0$, which is not the case. That there is at most
  one critical manifold of a given index is clear by an argument
  similar to the one in the proof of \reflemma{nonor-claim}.
\end{proof}

\section{Orientability of Negative Bundles}\label{allorient}

We want to use the fact that $\E$ is perfect to conclude that all
negative bundles over the three dimensional critical manifolds are
indeed orientable. Notice that this is not a circular argument, since
the proof of perfectness does not use orientability.

\begin{prop}\label{et-dim-udeluk}
  There do not exist any one dimensional critical manifolds.
\end{prop}

\begin{proof}
  This is similar to the argument in the proof of
  \reflemma{nonor-claim}. By $S^1$-equivariant perfectness of $\E$
  there can be at most two circles of index $2j$. Let $2k$, $k \geq 2$
  be the minimal index of a one dimensional critical manifold.  This
  critical manifold contributes a $\Q^2$ in degree $2k$. We must then
  fill the $2k -3$ gaps from $H_{S^1}^2(\Lambda S^3, S^3; \Q)$ to
  $H_{S^1}^{2k-2}(\Lambda S^3, S^3; \Q)$ with contributions that come
  in pairs or quartets, clearly impossible. If $k=1$ we get a
  contradiction since $H_{S^1}^2(\Lambda S^3, S^3; \Q) =\Q$
\end{proof}

\begin{cor}\label{oriented}
  The negative bundles over the three dimensional critical manifolds
  are ($S^1$-equivariantly) orientable.
\end{cor}

\begin{proof}
  First note that if the dimension of a critical manifold $N$ is
  three, the action of $S^1/\Z_n$ on $N$ is free, since there are no
  one dimensional critical manifolds.

  Hence the $S^1$-Borel construction $N \times_{S^1}ES^1 $ is homotopy
  equivalent to $N/(S^1/\Z_n)$ and by the proof of \refprop{eqcofix}
  $N/(S^1/\Z_n)$ is homotopy equivalent to $S^2$. Since $S^2$ is
  simply connected the bundle is orientable. We conclude that the
  negative bundles are also oriented as ordinary vector bundles.
\end{proof}

It was shown in the previous section that the negative bundles over
the one and five dimensional critical manifolds are orientable, so
this corollary finishes the proof that all negative bundles are
orientable.

\section{Topology of the Three Dimensional Critical Manifolds}

By calculating the Euler class of the bundle $S^1 \to N \to N/S^1 =
S^2$ we will deduce that the three dimensional critical manifolds are
integral cohomology three spheres or lens spaces.

\begin{thm}\label{cohomindex2}
  Assume that $N$ is a three dimensional critical manifold. The
  quotient $N/S^1 = S^2$ is endowed with a symplectic structure which
  corresponds to the Euler class of the $S^1$-bundle $S^1 \to N \to
  N/S^1$; in particular the Euler class is nonzero and $N$ has the
  integral cohomology of either the three sphere or a lens space.
\end{thm}

\begin{proof} 
  Note that by \refcor{conn} $N$ is connected. Parts of the proof rely
  on an argument in \cite{Be;mgc}; see \cite[Definition 1.23,
  Proposition 2.11 and 2.16]{Be;mgc}. We will describe some extra
  structure on the fixed point set $N$. In \cite[Chapter 2]{Be;mgc}
  the author considers a manifold all of whose geodesics are closed
  with the same least period $2\pi$. In that case the action of $S^1$
  on $T^1S^3$ is free and one gets a principal $S^1$ bundle $p\colon
  T^1S^3 \to T^1S^3/S^1$. There is a canonical connection $\alpha \in
  H^1_{\text{dR}}(T^1S^3; \mathcal{L}(S^1))$ on $T^1S^3$ constructed
  as follows: Let $p_{TS^3}\colon TTS^3 \to TS^3$ be the projection
  and $T_{p_{S^3}}\colon TTS^3 \to TS^3$ be the tangent map. For $X
  \in TTS^3$ we define $\tilde \alpha(X) = g(T_{p_{S^3}}(X),
  p_{TS^3}(X))$. This form is the pullback to the tangent bundle of
  the canonical one form on the cotangent bundle. Define a horizontal
  distribution on $TT^1S^3$ by $Q_u = \{X \in T_uT^1S^3 \mid \alpha
  (X) =0 \}$. Then $T_uT^1S^3 = \R Z \oplus Q_u$, where $Z$ is the
  geodesic vector field on $TT^1S^3$, and the corresponding connection
  form is $\tilde \alpha$ restricted to $T^1S^3$ which we denote by
  $\alpha$. The Lie algebra of $S^1$ is abelian and $d\alpha$ is
  horizontal, so the curvature form of $\alpha$ is $d\alpha$. Since
  the two form $d\alpha$ is invariant under the action of $S^1$, it is
  basic. The form $d \alpha$ is also the restriction to the unit
  tangent bundle of the pullback to the tangent bundle of the
  canonical two form on the cotangent bundle. Hence $d\alpha$ is
  nonzero on the complement of $Z$ on every $T_uT^1S^3$, $u \in
  T^1S^3$. By Chern-Weil Theory there exists $\omega \in
  H_{\text{dR}}^2(T^1S^3/S^1)$ such that $p^*(\omega) =d\alpha$. By
  \cite[Theorem 5.1]{KN;2} this class is the Euler class of the
  bundle. By \cite[Proposition 2.11]{Be;mgc} the class $\omega$ is a
  symplectic form on $T^1S^3/S^1$ ($\omega$ being degenerate would
  mean that $d \alpha =0$ on a nonempty subset of the complement of
  $Z$, which is not the case). Hence the Euler class is nonzero.

  We now carry the argument over to the three dimensional critical
  manifolds. By \refprop{et-dim-udeluk} we know that the action of
  $S^1$ on $N$ is free and that $N/S^1$ is homeomorphic to
  $S^2$. Consider the principal bundle $S^1 \to N \to S^2$, with
  projection $q$. We want to conclude that the Euler class of this
  bundle is nonzero. First note that the geodesic vector field is
  tangent to the fixed point set. Consider the restriction of $\alpha$
  to $TN$ and define as above a distribution on $TN$ by $\tilde Q_u =
  \{X \in T_uN \mid \alpha(X) =0\}$. Then we have $TN = \R Z \oplus
  \tilde Q_u$, since $Z$ is a vector field on $TN$ with $\alpha (Z) =
  1$; see \cite[1.57]{Be;mgc}. The restriction of $\alpha$ to $TN$ is
  invariant under the action of $S^1$ given by the geodesic flow,
  since by \cite[1.56]{Be;mgc} $L_Z \alpha =0$, so the distribution
  $\tilde Q$ is invariant under the action of $S^1$. Since the map $u
  \mapsto \tilde Q_u$ is clearly smooth, $\tilde Q$ is a horizontal
  distribution and since $\alpha(Z) =1$, $\alpha$ is the connection of
  the distribution. Similarly to the above, we see that since the Lie
  algebra of $S^1$ is abelian the curvature of the bundle is $d\alpha$
  ($d\alpha$ is horizontal since $d\alpha(Z, - ) =0$ by
  \cite[1.56]{Be;mgc}). As $d\alpha$ is horizontal and invariant under
  the action of $S^1$ ($L_Zd\alpha =0$ by \cite[1.57]{Be;mgc}), it is
  basic. Hence we can find a form $\omega \in H^2_{\text{dR}}(S^2)$
  such that $q^*(\omega) = d\alpha$. Thus, to see that the quotient is
  symplectic, it suffices to show that the curvature is nonzero on
  $\tilde Q$. This follows from the following general statement about
  symplectic reduction. The proof of this statement is that, similarly
  to the proof of \reflemma{fixorient}, the $+1$ eigenspace of the
  Poincar\'{e} map is a symplectic subspace.

  \begin{lemma}
    Let $(V^4, \tau)$ be a symplectic vector space and let $\Z_n$ act
    on $V$ by linear symplectic transformations. Assume that $\dim
    \Fix(\Z_n) =2$. Then $\Fix(\Z_n)$ is a symplectic subspace.
  \end{lemma}

  Consider the tangent space to $T^1S^3$ at a point $u \in
  T^1S^3$. This splits as a direct sum $Z \oplus Q_u$ and $d \alpha
  \neq 0$ on $Q_u$ and on the four dimensional subspace the form $d
  \alpha$ is a symplectic two form and the differential of the
  geodesic flow acts by symplectic linear transformations. The tangent
  space $T_uN$ also splits as the direct sum of $Z \oplus \tilde Q_u$,
  $\tilde Q_u \subset T_uN$. It follows from the above lemma that the
  form $d\alpha$ restricted to $T_uN$ is nonzero and hence, as above,
  that the form $\omega \in H^2_{\text{dR}}(S^2)$ makes the quotient
  into a symplectic manifold. Again by \cite[Theorem 5.1]{KN;2} the
  class $\omega$ is the Euler class of the $S^1$-bundle $S^1 \to N \to
  S^2$. Since the Euler class is nonzero an application of the Gysin
  sequence for the bundle $S^1 \to N \to S^2$ shows that $N$ is either
  an integral cohomology three sphere or lens space.
\end{proof}

\begin{cor}\label{reven}
  There are no three dimensional critical manifolds with the integral
  cohomology of $S^3$ or a lens space $S^3/\Z_r$, $r$ odd.
\end{cor}

\begin{proof}
  By \refcor{conn} there exists at most one connected, critical
  manifold, $N$, of a given index. Let $n$ be the multiplicity of a
  geodesic in $N$, i.e. $c$ has length $2\pi/n$. The action of $O(2)$
  leaves $N$ invariant and hence induces an action on $N$ which is
  effectively free since there are no one dimensional critical
  manifolds. Let $\Z_n \subseteq S^1$ be the ineffective kernel of the
  action. By identifying $S^1/\Z_n$ with $S^1$ we see that $O(2)/\Z_n
  \cong O(2)$. The group $O(2)/\Z_n$ acts freely on $N$, so in
  particular $\Z_2 \times \Z_2 \subseteq O(2)/\Z_n$ acts freely on
  $N$. By a a theorem of Smith \cite[Theorem 8.1]{Br;compgr}
  $\Z_2$-cohomology spheres do not support a free action of $\Z_2
  \times \Z_2$, so since $S^3/\Z_r$, $r$ odd, are $\Z_2$-cohomology
  spheres we have finished the proof.
\end{proof}

\section{Geometric Consequences}\label{conseq}

In this section we note some geometric consequences of the above
results. The first consequence is that for $g$ a metric on $S^3$ all
of whose geodesics are closed, the geodesics have the same least
period if and only if the energy function is perfect for ordinary
cohomology.

Recall that we assume the existence of exceptional shorter geodesics
of period $2\pi/n$, $n> 1$. These exceptional geodesics must lie in a
three dimensional critical manifold, since we have seen that there
exist no one dimensional critical manifolds, and if the geodesics of
period $2\pi/n$ lie in a five dimensional critical manifold, all
geodesics are closed of period $2\pi/n$, contrary to our
assumption. 

We consider the unique critical manifold of index two, which we denote
by $N$. To prove the conjecture there are then two cases to rule out:
$N$ is three and five dimensional. The five dimensional case is easy
since the shorter geodesics also have index two, contradicting the
fact that there is at most one critical manifold of each given index.

For $N$ three dimensional and $c \in N$ we have two cases to consider:
$\ind(c^2) = 4$ and $\ind(c^2) = 6$. If we assume that the sectional
curvature $K$ is pinched $a/4 \leq K \leq a$, it follows from
\cite{BTZ;CLOPOS} and \cite{Tsuka;clo} that the critical manifolds
contributing the classes in degree $2$, $4$, $5$ and $7$ to
$H^*(\Lambda S^3; \Z)$ only contain prime closed geodesics. Hence we
get a contradiction in the case where $\ind(c^2) = 4$. We have not yet
been able to handle the case where $\ind(c^2) = 6$.

\newpage

\bibliography{litt}

\end{document}